\documentclass[10pt,twoside]{article}
\usepackage{graphicx}
\usepackage{Latex-document}

\newcommand{\pictype}{}             

\newlength{\dddotwidth}
\newlength{\dddotwidthb}
\newlength{\dddotheight}
\newcommand{\dddot}[1]{{
  \settowidth{\dddotwidth}{\ensuremath{#1}}
  \settowidth{\dddotwidthb}{$\dot{\,}\dot{\,}\dot{\,}$}
  \addtolength{\dddotwidth}{\dddotwidthb}
  \settoheight{\dddotheight}{\ensuremath{#1}}
  \ensuremath{
    \dot{\rule{0pt}{\dddotheight}\,}
    \dot{\rule{0pt}{\dddotheight}\,}
    \dot{\rule{0pt}{\dddotheight}\,}
    \hspace{-0.5\dddotwidth}#1}
  }}

\newcommand{\setR}{\mathsf{I\hspace{-0.11em}R}}

\renewcommand{\thesection}{\arabic{section}}

\def\addsec{}

\def\Section#1{\section{\hskip -1em . \hskip 0.8em #1}}

\title{\bf Bifurcations without Parameters:\vskip -2mm Some ODE and PDE Examples\vskip 6mm}

\author{{\bf Bernold Fiedler}\thanks{Free University Berlin, Institute of Mathematics I, Arnimallee 2-6, D-14195 Berlin,
Germany. E-mail: fiedler@math.fu-berlin.de} \quad Stefan
Liebscher\thanks{Free University Berlin, Institute of Mathematics
I, Arnimallee 2-6, D-14195 Berlin, Germany. E-mail:
liebsch@math.fu-berlin.de, www.math.fu-berlin.de/\~{ }Dynamik/}
\vspace*{-0.5cm}}

\date{\vspace{-8mm}}

\begin{document}

\maketitle

\thispagestyle{first} \setcounter{page}{305}

\begin{abstract}

\vskip 3mm

Standard bifurcation theory is concerned with families of vector
fields
$$ \dot x = f(x,\lambda) $$
$x \in  {\bf R}^n$, involving one or several constant real
parameters $\lambda$. Viewed as a differential equation for the
pair $(x,\lambda)$, we observe a foliation of the total phase
space by constant $\lambda$. Frequently, the presence of a trivial
stationary solution $x=0$ is also imposed:
$$0 = f(0,\lambda) .$$

Bifurcation without parameters, in contrast, discards the
foliation by a constant parameter $\lambda$. Instead, we consider
systems
 \begin{eqnarray*}
        \dot x &=& f(x,y) \\
        \dot y &=& g(x,y)
        \end{eqnarray*}
Standard bifurcation theory then corresponds to the special case
$y=\lambda,\ g=0$. To preserve only the trivial solution $x=0$,
instead, we only require
$$ 0 = f(0,y) = g(0,y) $$
for all $y$. A rich dynamic phenomenology arises, when normal
hyperbolicity of the trivial stationary manifold $x=0$ fails, due
to zero or purely imaginary eigenvalues of the Jacobian
$f_x(0,y)$.

Specifically, we address the cases of failure of normal
hyperbolicity due to a simple eigenvalue zero, a simple purely
imaginary pair (Hopf bifurcation without parameters), a double
eigenvalue zero (Takens-Bogdanov bifurcation without parameters),
and due to a double eigenvalue zero with additional time reversal
symmetries.

We illustrate consequences of our results with examples from
ordinary and partial differential equations arising in systems of
coupled oscillators, in the analysis and numerics of hyperbolic
conservation and balance laws, and in the fluid dynamics of plane
Kolmogorov flows.

The results are joint work with Andrei Afendikov, James C.
Alexander, and Stefan Liebscher. For references see
\begin{center} www.math.fu-berlin.de/$\sim$Dynamik/
\end{center}

\vskip 4.5mm


\end{abstract}

\vskip 12mm

\Section{Applied motivation} \label{sec1} \setzero \addsec

\vskip-5mm \hspace{5mm}

In this article we sketch and illustrate some elements of the nonlinear
dynamics near equilibrium manifolds. Denoting the equilibrium manifold by
$x=0$, in local coordinates $(x,y)\in\setR^n\times \setR^k$, we consider
systems
\begin{equation}\begin{array}{rcl}
   \dot{x} &=& f(x,y)
\\ \dot{y} &=& g(x,y)
\label{1.1}\end{array}\end{equation}
and assume
\begin{equation}\begin{array}{rclcl}
   f(0,y)&=&g(0,y)&=&0,
\label{1.2}\end{array}\end{equation}
for all $y$. For simplicity we will only address the cases $k=1$ of lines of
equilibria, and $k=2$ of equilibrium planes. Sufficient smoothness of $f,g$ is
assumed. The occurrence of equilibrium manifolds is infinitely degenerate, of
course, in the space of all vector fields $(f,g)$ -- quite like many
mathematical structures are: equivariance under symmetry groups, conservation
laws, integrability, symplectic structures, and many others. The special case
\begin{equation}\begin{array}{rcl}
   g &\equiv& 0
\label{1.3}\end{array}\end{equation}
in fact amounts to standard bifurcation theory, in the presence of a trivial
solution $x=0$; see for example \cite{ChowHale82-Bifurcation}.
Note that condition
(\ref{1.3}), which turns the $k$-dimensional variable $y$ into a preserved
constant parameter, is infinitely degenerate even in our present setting
(\ref{1.2}) of equilibrium manifolds. Due to the analogies of our results and
methods with bifurcation theory, we call our emerging theory {\em bifurcation
  without parameters}. This terminology emphasizes the intricate dynamics
which arises when normal hyperbolicity of the equilibrium manifold fails; see
the sections below.

To motivate assumption (\ref{1.2}), we present several examples. First,
consider an ``octahedral'' graph $\Gamma$ of $2(m+1)$ vertices $\{\pm
1,\ldots,\pm(m+1) \}$. The graph $\Gamma$ results from the complete graph by
eliminating the ``diagonal'' edges, which join the antipodal vertices $\pm j$,
for $j=1,\ldots,m+1$. For $m=1$ we obtain the square, for $m=2$ the
octahedron, and so on. Consider the system
\begin{equation}\begin{array}{rcl}
   \dot{u}_j &=& \displaystyle f_j(u_j, \sum_{k\neq \pm j} u_k)
\label{1.4}\end{array}\end{equation}
of oscillators $u_j\in\setR^{n'}$ on $\Gamma$, additively coupled along the
edges by $f_j$. We assume an antipodal oddness symmetry of the individual
oscillator dynamics
\begin{equation}\begin{array}{rcl}
   f_{-j} (-u_j,0) &=& -f_j(u_j, 0).
\label{1.5}\end{array}\end{equation}
As a consequence, the antipode space
\begin{equation}\begin{array}{rcl}
   \Sigma &:=& \{u=(u_j)_{j\in\Gamma};\;\;u_{-j}=-u_j \}
\label{1.6}\end{array}\end{equation}
is invariant under the flow (\ref{1.4}). Moreover, the flow on $\Sigma$
completely decouples into a direct product flow of the $m+1$ diagonally
antipodal, decoupled pairs
\begin{equation}\begin{array}{rcl}
   \dot{u}_{\pm j} &=& f_{\pm j} (u_{\pm j},0).
\label{1.7}\end{array}\end{equation}
For the square case $m=1$, this decoupling phenomenon was first observed in
\cite{AlexanderAuchmuty86-PhaseLockOsc}. For more examples see also
\cite{AlexanderFiedler89-Decoupling}.

An $m$-plane of equilibria arises from periodic solutions of the decoupled
system (\ref{1.7}). Assume (\ref{1.7}) possesses time periodic orbits
$u_j(t+\varphi_j)$ of equal period $T_j=2\pi$, for $j=1,\ldots,m+1$. Choose
arbitrary phases $\varphi_j\in S^1$ and let $u_{-j}(t):=
-u_j(t),\;\varphi_{-j}=\varphi_j$. Then
\begin{equation}\begin{array}{rclcl}
   u^\varphi (t) &:=& (u_j(t+\varphi_j))_{j\in \Gamma } &\in& \Sigma
\label{1.8}\end{array}\end{equation}
is a $2\pi$-periodic solution of (\ref{1.4}), (\ref{1.7}), for arbitrary
phases $\varphi \in T^{m+1}$. Eliminating one phase angle $\varphi_{m+1}$ by
passing to an associated Poincar\'e map, an $m$-dimensional manifold of fixed
points arises, parametrized by the remaining $m$ phase angles. Assuming, in
addition to the diagonal oddness symmetry (\ref{1.5}), equivariance of
(\ref{1.4}) with respect to an $S^1$-action, the Poincar\'e map can in fact be
obtained as a time-$2\pi$ map of an autonomous flow within the Poincar\'e
section. In suitable notation, $y=\{\varphi_1,\ldots,\varphi_m\}$, the fixed
point manifold then becomes an $m$-dimensional manifold of equilibria, as
presented in (\ref{1.1}), (\ref{1.2}) above. For more detailed discussions of
this example in the context of bifurcations without parameters see
\cite{Liebscher97-Diplom, FiedlerLiebscherAlexander98-HopfTheory,
  FiedlerLiebscher01-TakensBogdanov}.

As a second example of equilibrium manifolds we consider viscous profiles
$u=u((\xi -st)/\varepsilon)$ of systems of nonlinear hyperbolic conservation
laws and stiff balance laws
\begin{equation}\begin{array}{rcl}
   \partial_t u + \partial_{\xi} F(u) + \varepsilon^{-1} G(u) &=&
   \varepsilon \partial^2_\xi u.
\label{1.9}\end{array}\end{equation}
Viscous profiles then have to satisfy an $\varepsilon$-independent ODE system
\begin{equation}\begin{array}{rcl}
   \ddot{u} &=& (F'(u) - s \cdot \mathrm{id}) \dot{u} + G(u).
\label{1.10}\end{array}\end{equation}
Standard conservation laws, for example, require $G\equiv 0$. The presence of
$m$ conservation laws corresponds to nonlinearities $G$ with range in a
manifold of codimension $m$ in $u$-space. Typically, then, $G(u)=0$ describes
an equilibrium manifold of dimension $m$ of pairs $(u,\dot{u})=(u,0)$, in the
phase space of (\ref{1.10}). For an analysis of this example in the context of
bifurcations without parameters see \cite{FiedlerLiebscher98-HopfViscBalLaw,
Liebscher00-Dissertation}. For another
example, which relates binary oscillations in central difference
discretizations of hyperbolic balance laws with diagonal uncoupling of coupled
oscillators, see \cite{FiedlerLiebscherAlexander98-HopfBinOsc}.

We conclude our introductory excursion with a brief summary of some further
examples. In \cite{Farkas84-ZipBifurcation},
lines of equilibria have been observed for the
dynamics of models of competing populations. This included
a first partial analysis of failure of normal hyperbolicity.

A topologically very interesting example in compact three-dimensional
manifolds involves contact structures $\eta (\xi)$ (i.e., nonintegrable plane
fields and gradient vector fields $\dot{\xi} = -\nabla V(\xi) \in
\eta(\xi)$. See \cite{EtnyreGhrist99-GradientFlows} for an in-depth
analysis. Examples include mechanical systems with nonholonomic
constraints. Notably, level surfaces of regular values of the potential $V$
consist of tori. Under a nondegeneracy assumption, equilibria form embedded
circles, that is, possibly linked and nontrivial knots.

For a detailed study of plane Kolmogorov fluid flows in the presence of a line
of equilibria with a degeneracy of Takens-Bogdanov type and an additional
reversibility symmetry, see
\cite{AfendikovFiedlerLiebscher02-PlaneKolmogorovFlows}.

As a caveat we repeat that lines of equilibria, which are transverse to level
surfaces of preserved quantities $\lambda $ do {\em not} provide bifurcations,
without parameters. Rather, $\dot{y}=0$ for $y:=\lambda $ exhibits this
problem as belonging to standard bifurcation theory; see (\ref{1.3}).

\Section{Sample vector fields and resulting flows} \label{sec2} \setzero
\addsec

\vskip-5mm \hspace{5mm}

In this section we collect relevant example vector fields (\ref{1.1}),
(\ref{1.2}) with lines and planes of equilibria $x=0$; see
\cite{FiedlerLiebscherAlexander98-HopfTheory,
  FiedlerLiebscher01-TakensBogdanov,
  AfendikovFiedlerLiebscher02-PlaneKolmogorovFlows}
for further details. We illustrate and comment the resulting flows.

\begin{figure}
\setlength{\unitlength}{0.49\textwidth}
\begin{center}
\begin{picture}(1.0,1.0)(0,0)
  \put(0,0){\framebox(1.0,1.0)
   {
   \includegraphics[width=1.0\unitlength]{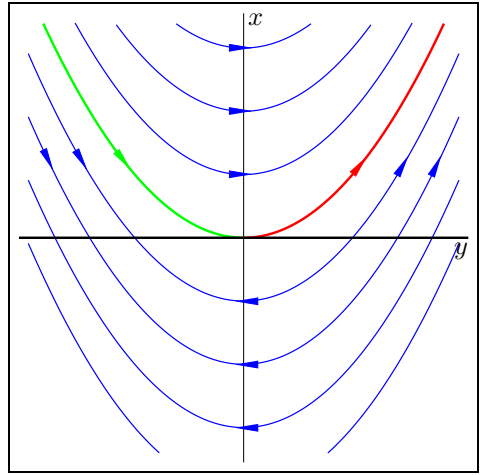}
   }}
  \put(0.98,0.49){\makebox(0,0)[rt]{$y$}}
  \put(0.51,0.98){\makebox(0,0)[lt]{$x$}}
\end{picture}
\begin{minipage}{10cm}
\caption{A line of equilibria ($y$-axis) with a nontrivial transverse
  eigenvalue zero.}
\label{fig2.1}
\end{minipage}
\end{center}
\end{figure}

Normally hyperbolic equilibrium manifolds admit a  transverse $C^0$-foliation
with hyperbolic linear flows in the leaves. See for example
\cite{Palmer84-ExpDichotHomPoint}, \cite{Shoshitaishvili75-Bif} and the ample
discussion in \cite{Aulbach84-EqManifolds}. As a first example, we therefore
consider
\begin{equation}\begin{array}{rcl}
   \dot{x} &=& x y,
\\ \dot{y} &=& x.
\label{2.1}\end{array}\end{equation}
Note the loss of normal hyperbolicity at $x=y=0$, due to a nontrivial
transverse eigenvalue zero of the linearization. Clearly $dx/dy=y$, and the
resulting flow lines are parabolas; see Figure \ref{fig2.1}. For comparison
with standard bifurcation theory, where $y=\lambda$, we draw the $y$-axis of
equilibria horizontally.

\begin{figure}
\setlength{\unitlength}{0.49\textwidth}
\begin{center}
\begin{picture}(1.0,1.0)(0,0)
  \put(0.02,0.98){\makebox(0,0)[lt]{(a)}}
  \put(0,0){\framebox(1.0,1.0)
   {
   \includegraphics[width=1.0\unitlength]{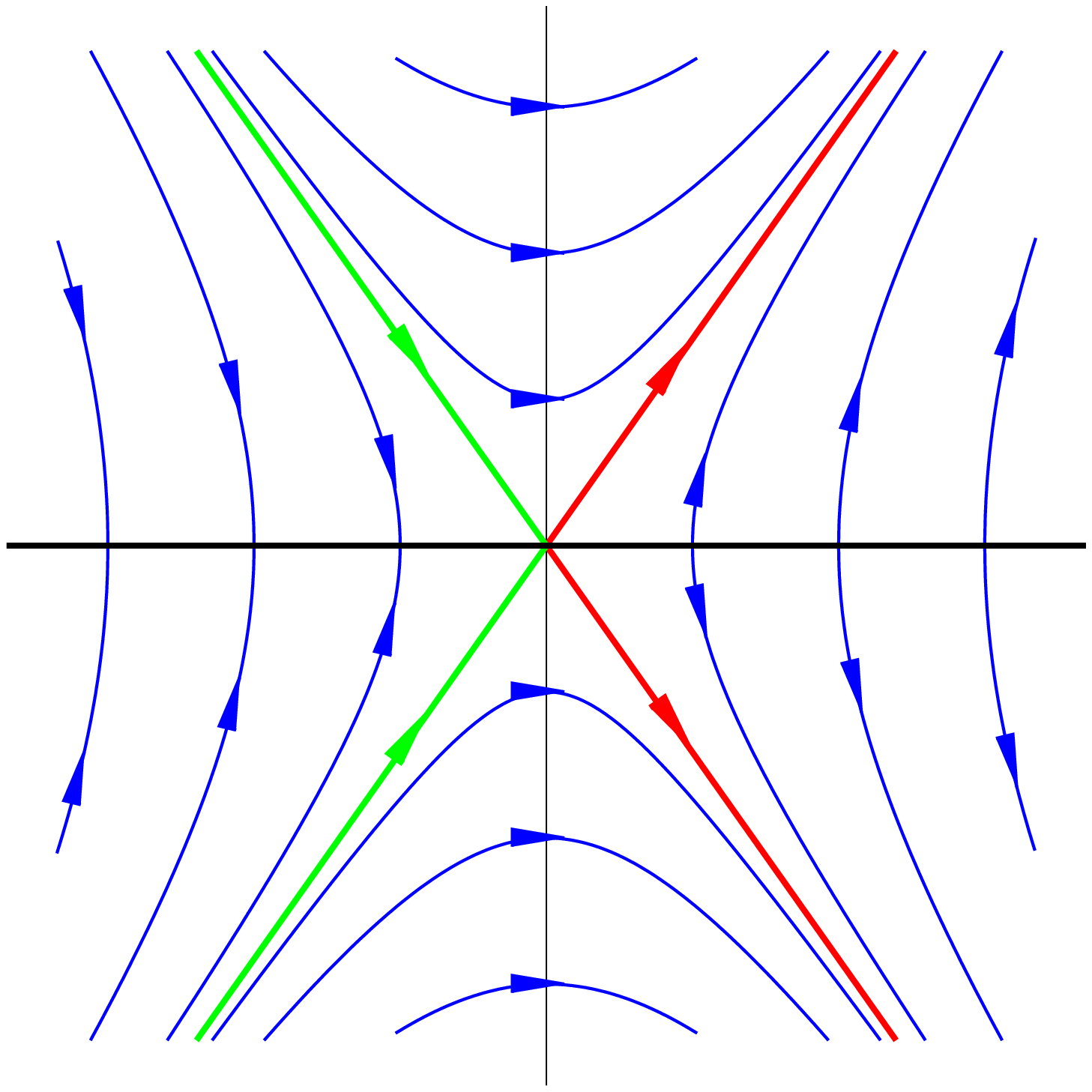}
   }}
  \put(0.97,0.48){\makebox(0,0)[rt]{$y$}}
  \put(0.51,0.97){\makebox(0,0)[lt]{$r$}}
\end{picture}
\begin{picture}(1.0,1.0)(0,0)
  \put(0,0){\makebox(1.0,1.0)
   {
   \includegraphics[width=1.0\unitlength]{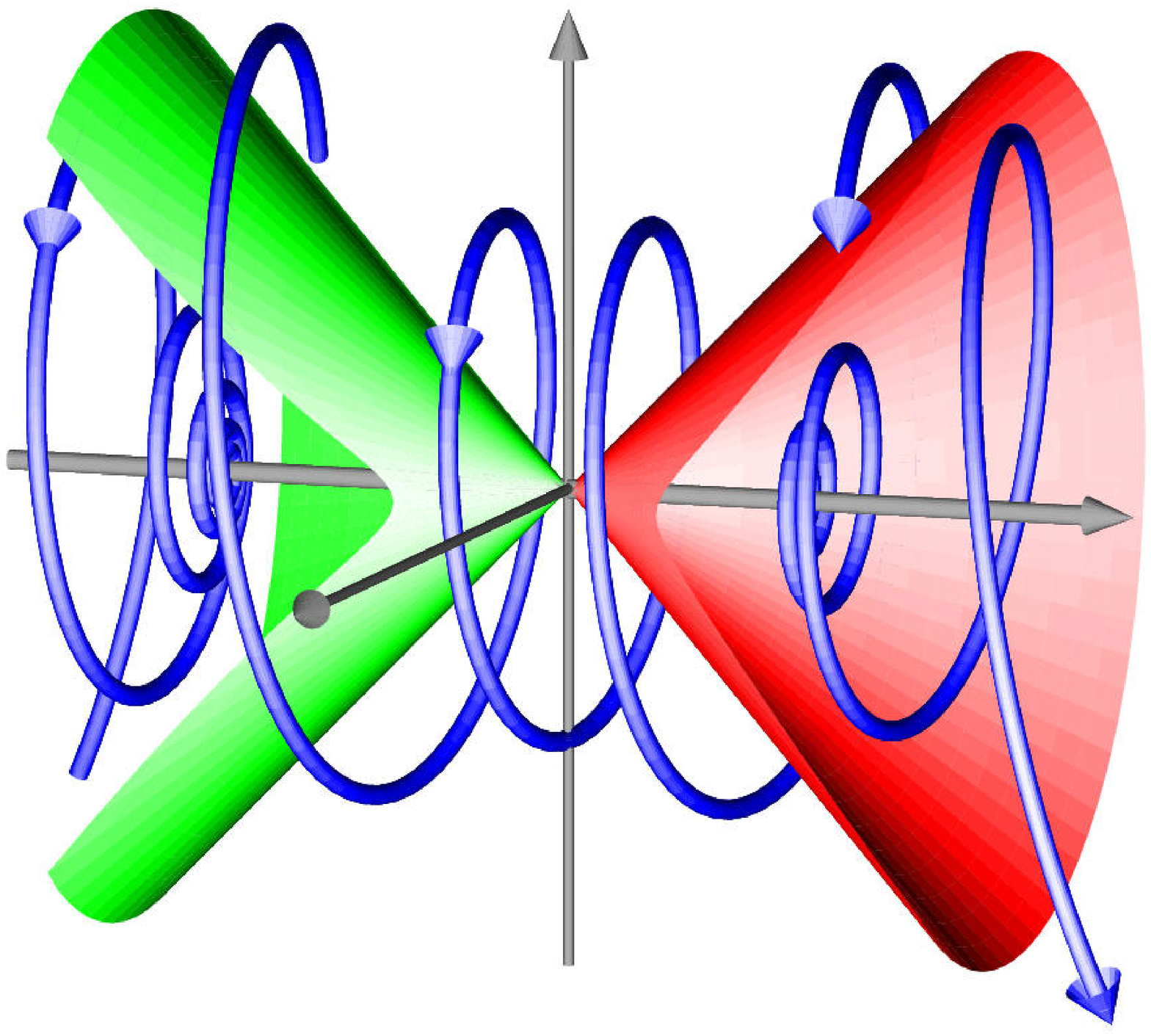\pictype}
   }}
  \put(0,0){\framebox(1.0,1.0){~}}
  \put(0.97,0.45){\makebox(0,0)[rt]{$y$}}
  \put(0.52,0.90){\makebox(0,0)[lt]{$x_1$}}
  \put(0.29,0.38){\makebox(0,0)[lt]{$x_2$}}
\end{picture}
\\[1ex]
\begin{picture}(1.0,1.0)(0,0)
  \put(0.02,0.98){\makebox(0,0)[lt]{(b)}}
  \put(0,0){\framebox(1.0,1.0)
   {
   \includegraphics[width=1.0\unitlength]{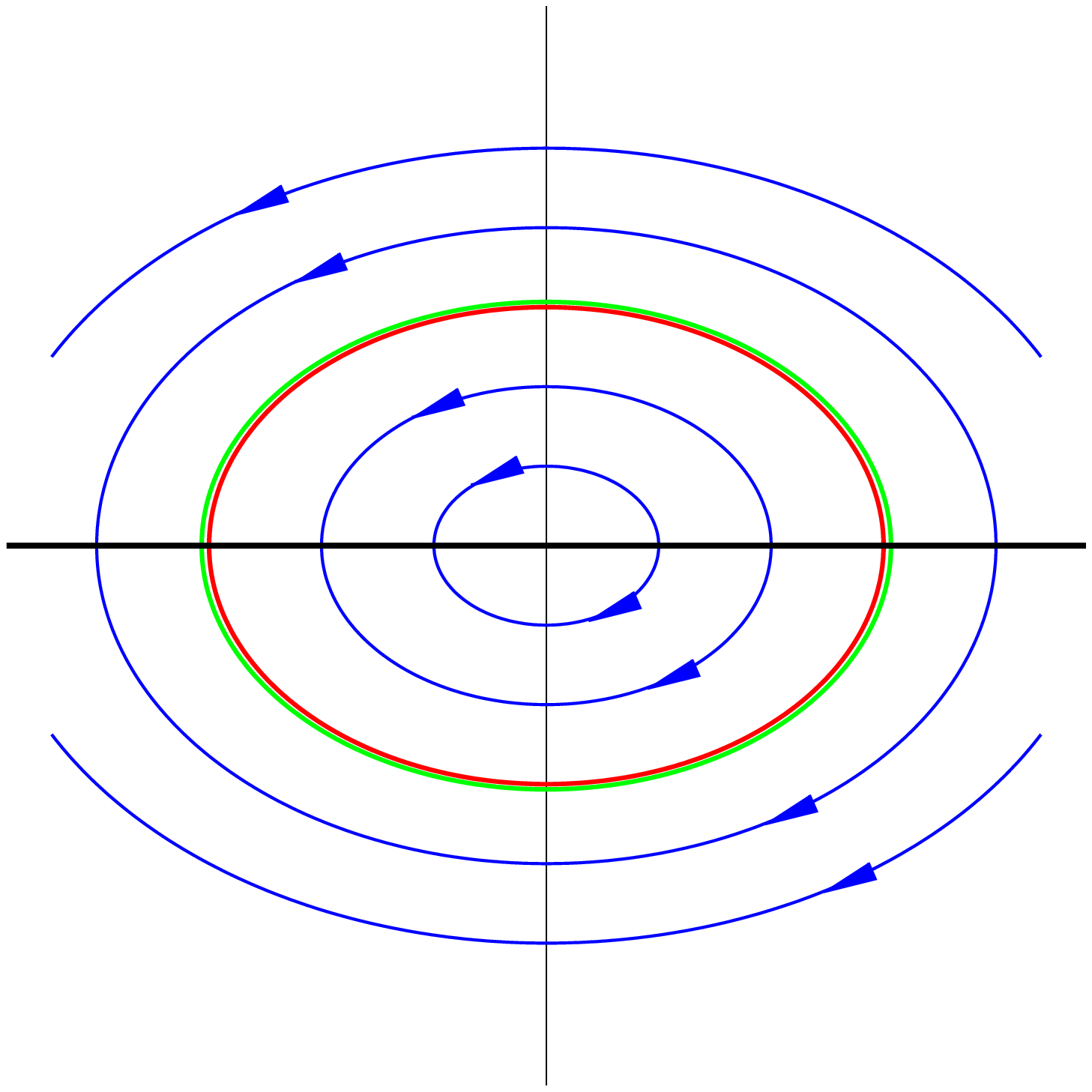}
   }}
  \put(0.97,0.48){\makebox(0,0)[rt]{$y$}}
  \put(0.51,0.97){\makebox(0,0)[lt]{$r$}}
\end{picture}
\begin{picture}(1.0,1.0)(0,0)
  \put(0,0){\makebox(1.0,1.0)
   {
   \includegraphics[width=1.0\unitlength]{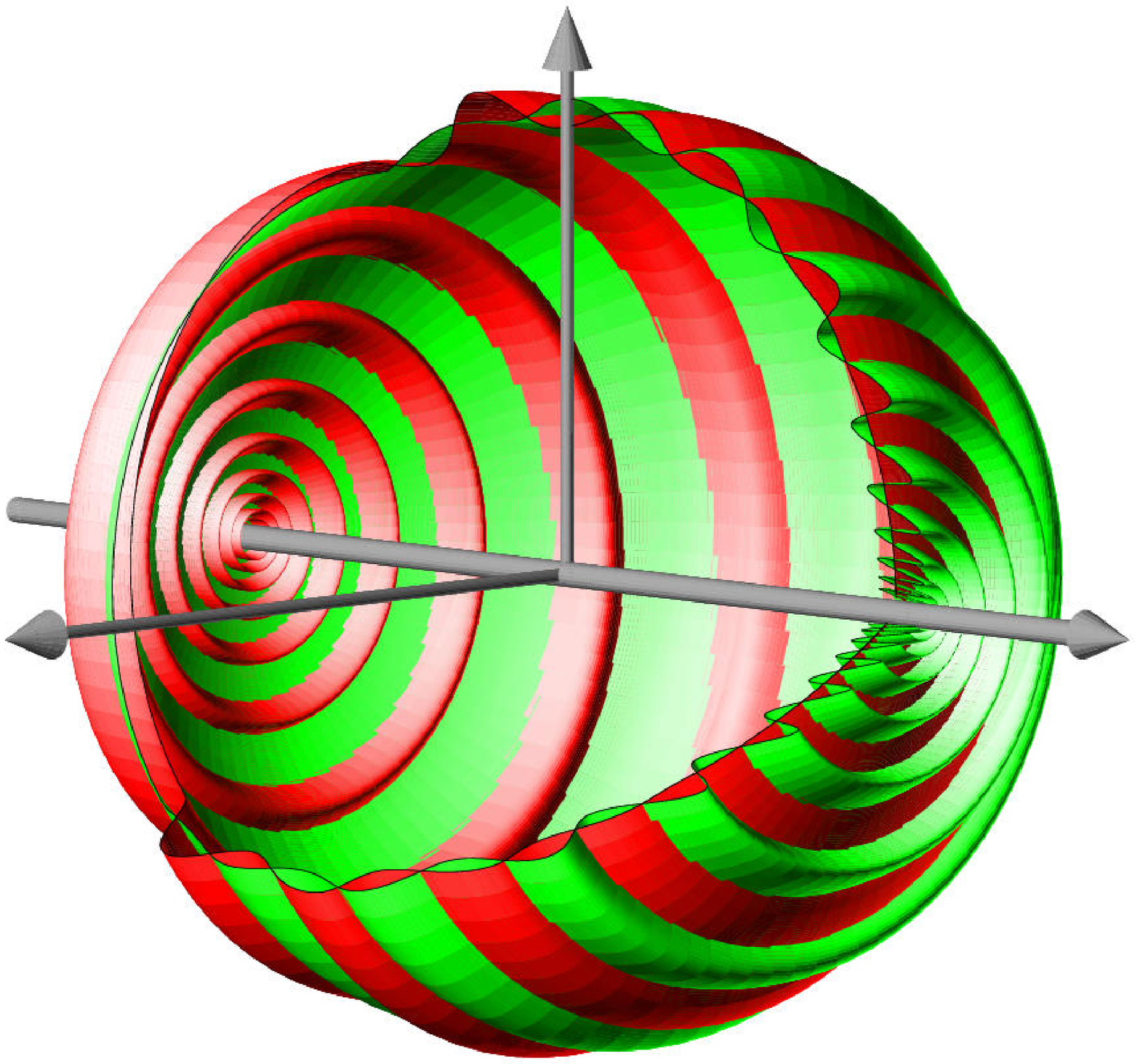\pictype}
   }}
  \put(0,0){\framebox(1.0,1.0){~}}
  \put(0.97,0.40){\makebox(0,0)[rt]{$y$}}
  \put(0.53,0.96){\makebox(0,0)[lt]{$x_1$}}
  \put(0.03,0.40){\makebox(0,0)[lt]{$x_2$}}
\end{picture}
\begin{minipage}{10cm}
\caption{Lines of equilibria ($y$-axis) with imaginary eigenvalues: Hopf
   bifurcation without parameters. Case (a) hyperbolic; case (b)
   elliptic. Red: strong unstable manifolds; green: stable manifolds.}
\label{fig2.2}
\end{minipage}
\end{center}
\end{figure}

As a second example, consider
\begin{equation}\begin{array}{rcl}
   \dot{x} &=& x y,
\\ \dot{y} &=& \pm x^2.
\label{2.2}\end{array}\end{equation}
Again, a transverse zero eigenvalue occurs -- this time with an additional
reflection symmetry $y \mapsto -y$. Dividing by the Euler multiplier $x$, the
reflection becomes a time reversibility. See the left parts of Figure
\ref{fig2.2} for the resulting flows. Note the resulting integrable, harmonic
oscillator case which originates from the elliptic sign $\dot{y}=-x^2$.

As a third example, we consider $x=(x_1,x_2)\in\setR^2$, $y\in \setR$ with a
line $x=0$ of equilibria and a purely imaginary nonzero eigenvalue $i\omega$
at $x=0$. Normal-form theory, for example as in
\cite{Vanderbowhede89-CentManifolds}, then generates an additional
$S^1$-symmetry by the action of $\exp(i\omega t)$ in the $x$-eigenspace. This
equivariance can be achieved, successively, up to Taylor expansions of any
finite order. In polar coordinates $(r,\varphi)$ for $x$, an example of
leading order terms is given by
\begin{equation}\begin{array}{rcl}
   \dot{r} &=& r y,
\\ \dot{y} &=& \pm r^2,
\\ \dot{\varphi} &=& \omega.
\label{2.3}\end{array}\end{equation}
Since the first two equations in (\ref{2.3}) coincide with (\ref{2.2}), the
dynamics is then obtained by simply rotating the left parts of Figure
\ref{fig2.2} around the $y$-axis at speed $w$. The right parts of Figure
\ref{fig2.2} provide three-dimensional views of the effects of higher-order
terms which do not respect the $S^1$-symmetry of the normal forms. In the
elliptic case (b), all nonstationary orbits are heteroclinic from the unstable
foci, at $y>0$, to stable foci at $y<0$. The two-dimensional respective strong
stable and unstable manifolds will split, generically, to permit transverse
intersections.

\begin{figure}
\begin{center}
\setlength{\unitlength}{0.69\textwidth}
\begin{picture}(1.44,2.0)(-0.72,-1.0)
  \put(0,0){\rotatebox[origin=c]{90}{\makebox(0,0){
  \begin{minipage}{2.03\unitlength}
\begin{picture}(1.0,0.8)(0,0)
  \put(0,0){\framebox(1.0,0.8)
   {
   \includegraphics[width=1.0\unitlength]{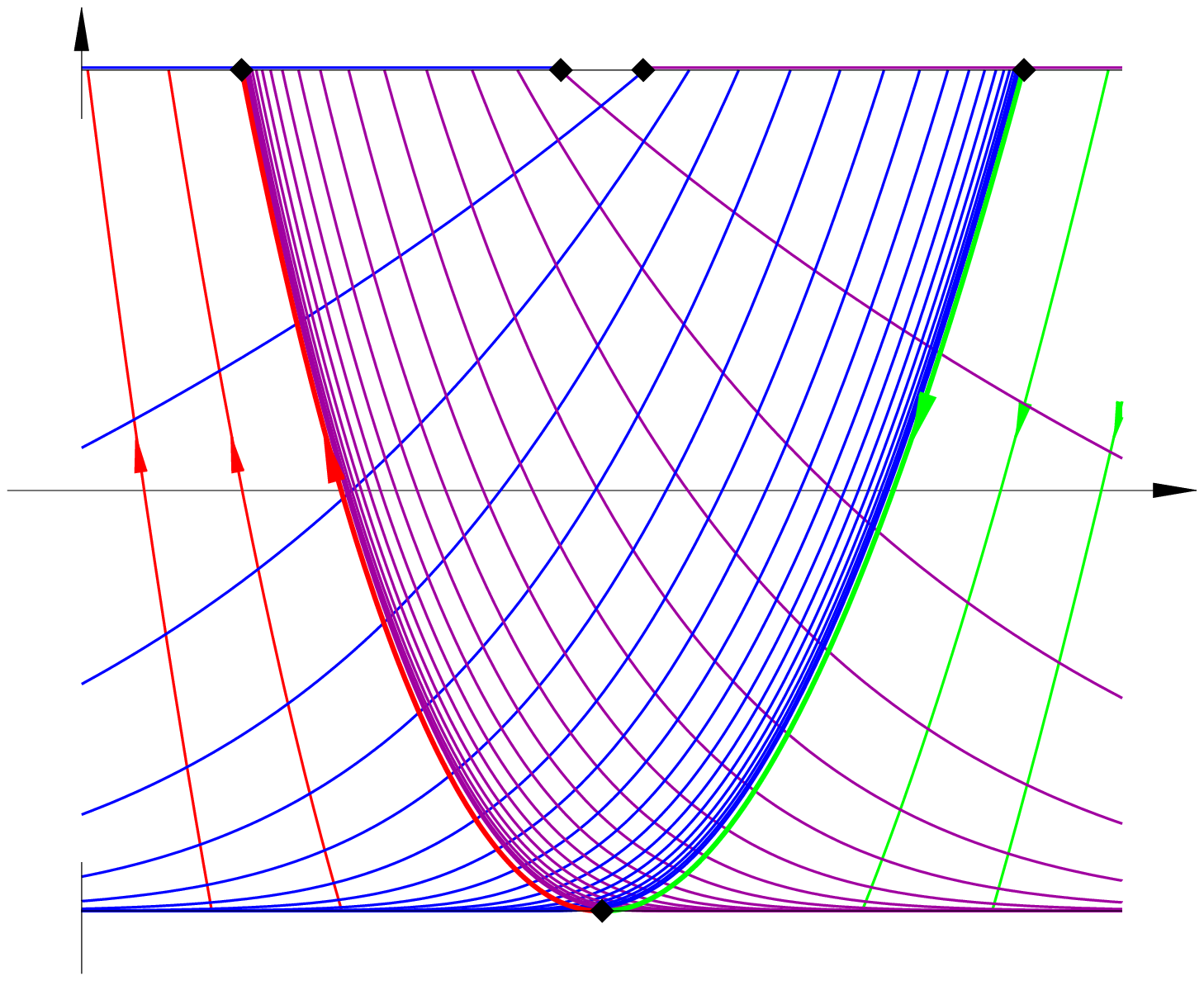}
   }}
  \put(0.5  ,0.055){\makebox(0,0)[t]{hyp.~Hopf\strut}}
  \put(0.495,0.790){\makebox(0,0)[t]{\small1-het}}
  \put(0.075,0.064){\makebox(0,0)[r]{\small$-\frac{\sqrt{8}}{3}$}}
  \put(0.075,0.745){\makebox(0,0)[r]{\small$\frac{\sqrt{8}}{3}$}}
  \put(0.09 ,0.792){\makebox(0,0)[lt]{$\tilde{H}$\strut}}
  \put(0.99 ,0.4  ){\makebox(0,0)[rt]{$\tau$\strut}}
  \put(0.5  ,0.6  ){\makebox(0,0)[t]{
    \shortstack{infinite\\heteroclinic\strut\\swarm}}}
\end{picture}
\hfill
\begin{picture}(1.0,0.8)(0,0)
  \put(0,0){\framebox(1.0,0.8)
   {
   \includegraphics[width=1.0\unitlength]{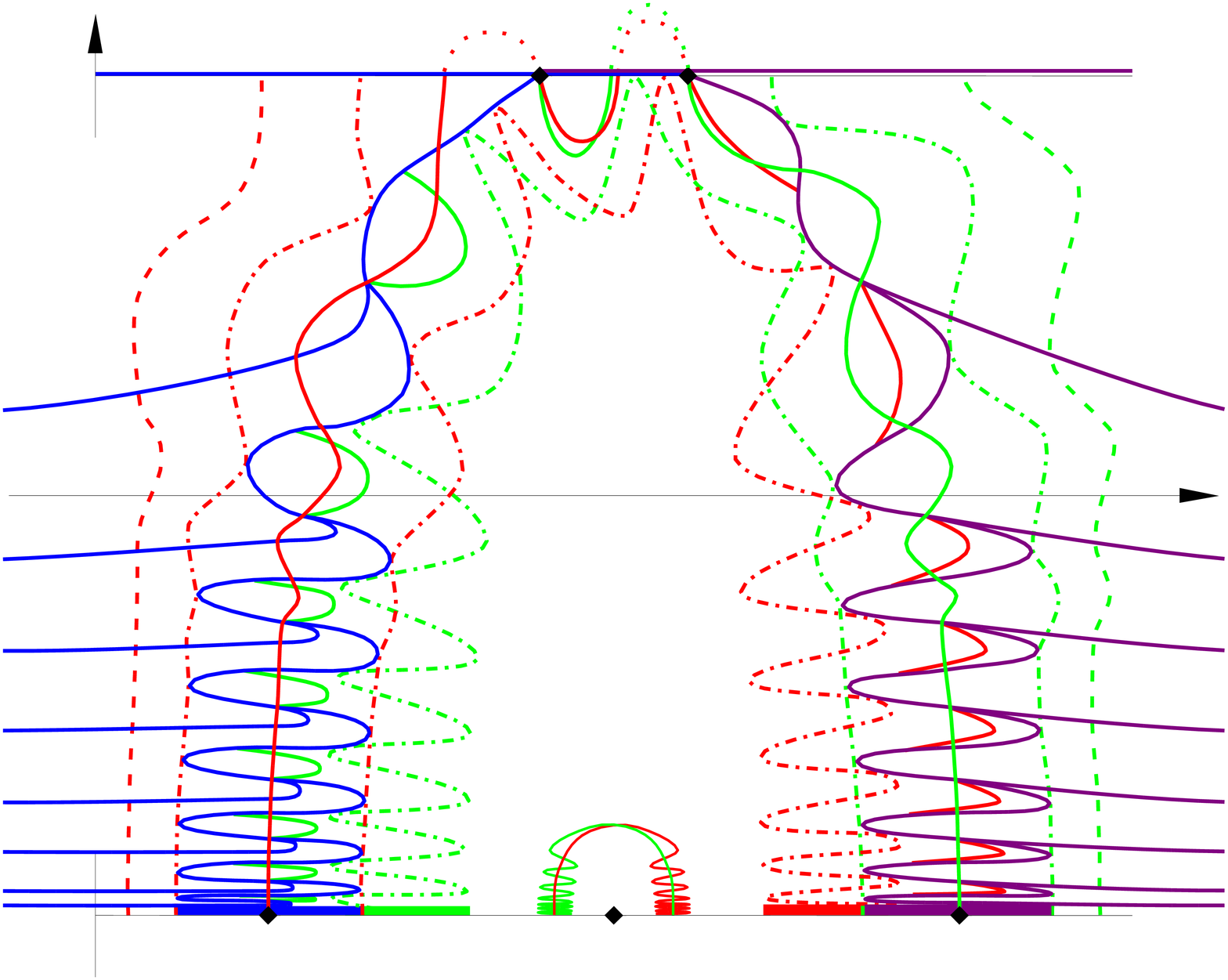}
   }}
  \put(0.5  ,0.055){\makebox(0,0)[t]{ell.~Hopf\strut}}
  \put(0.495,0.790){\makebox(0,0)[t]{\small1-het}}
  \put(0.075,0.045){\makebox(0,0)[r]{\small$-\frac{\sqrt{8}}{3}$}}
  \put(0.075,0.745){\makebox(0,0)[r]{\small$\frac{\sqrt{8}}{3}$}}
  \put(0.09 ,0.792){\makebox(0,0)[lt]{$\tilde{H}$\strut}}
  \put(0.99 ,0.4  ){\makebox(0,0)[rb]{$\tau$\strut}}
\end{picture}
\\[1ex]
\begin{picture}(1.0,0.5)(0,0)
  \put(0,0){\framebox(1.0,0.5)
   {
   \includegraphics[width=1.0\unitlength]{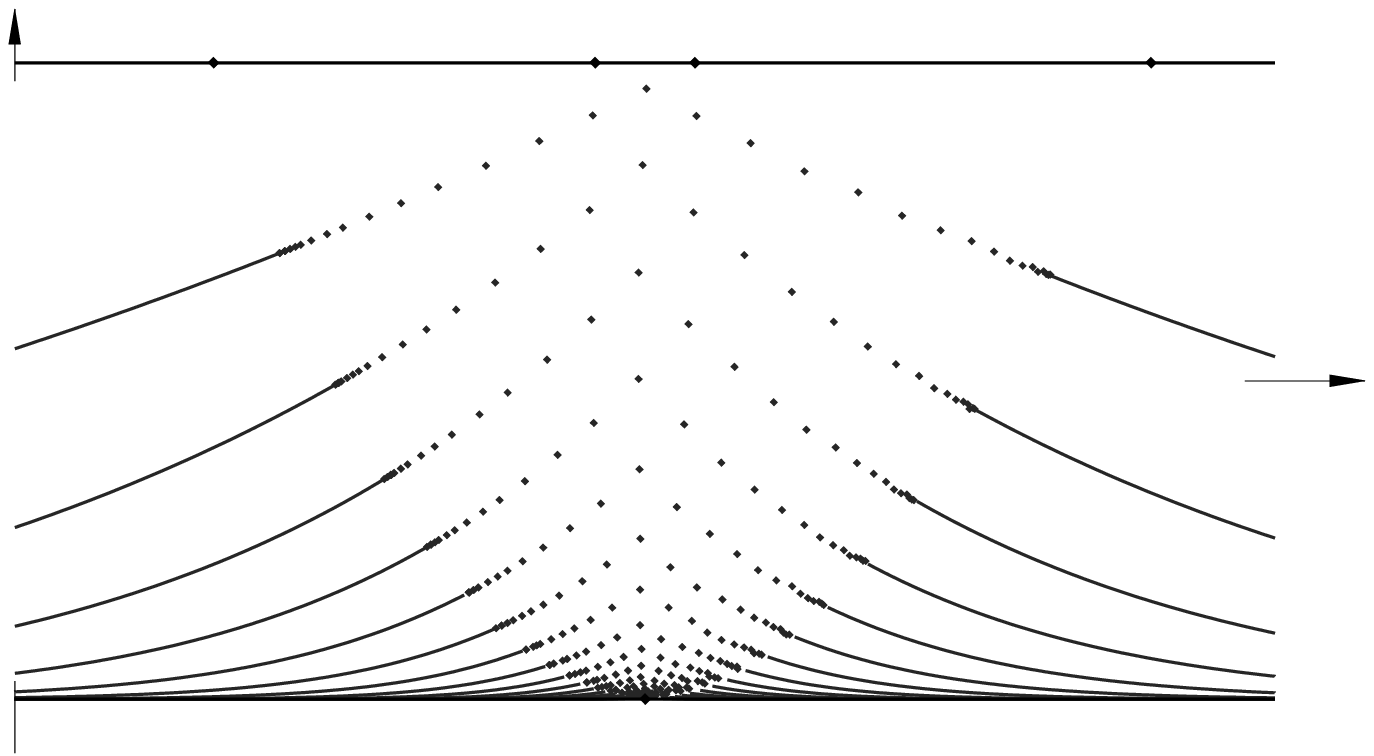}
   }}
  \put(0.07,0.495){\makebox(0,0)[rt]{$\tilde{H}$\strut}}
  \put(0.99,0.25 ){\makebox(0,0)[rt]{$\tau$\strut}}
\end{picture}
\hfill
\begin{picture}(1.0,0.5)(0,0)
  \put(0,0){\framebox(1.0,0.5)
   {
   \includegraphics[width=1.0\unitlength]{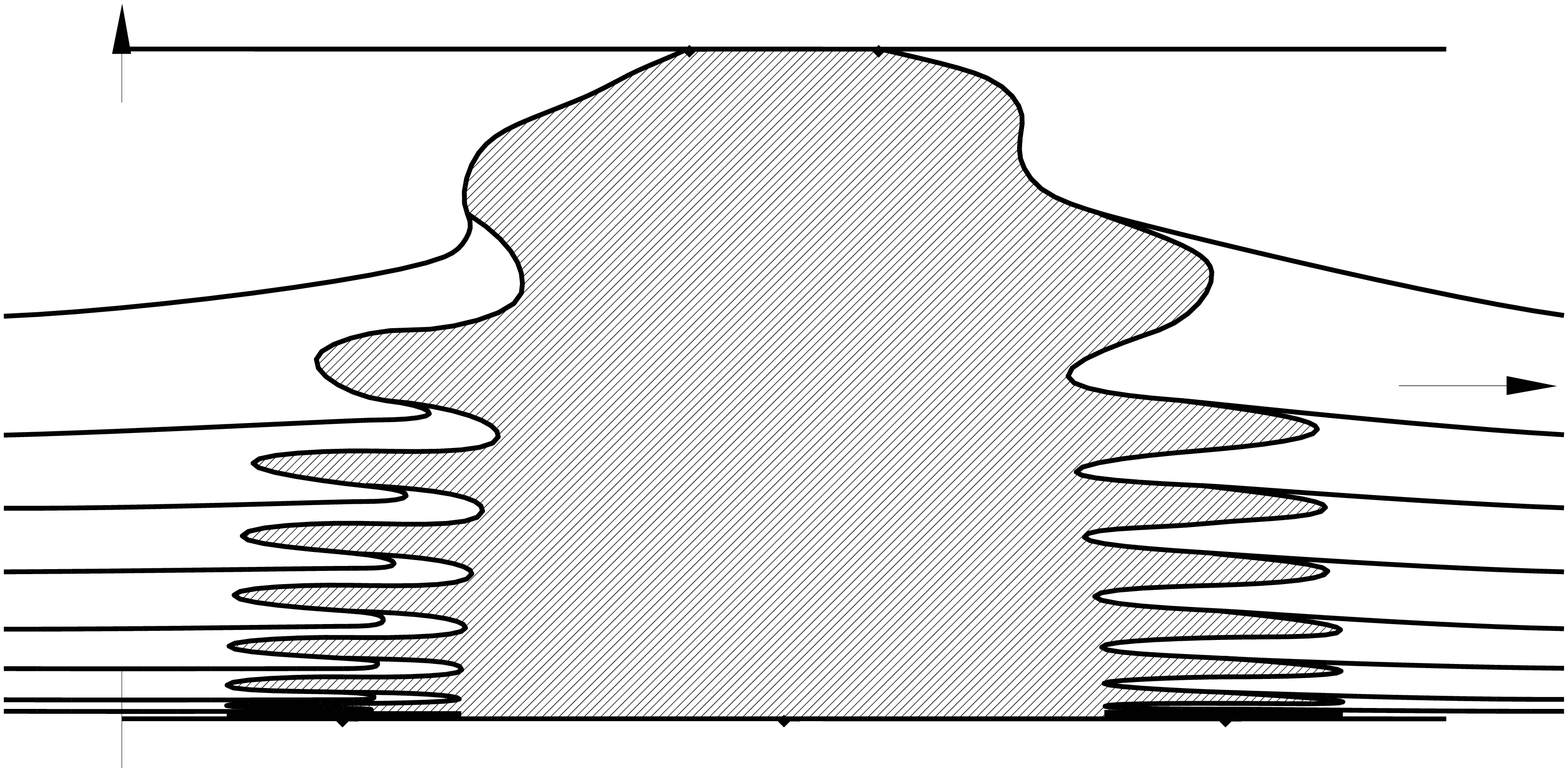}
   }}
  \put(0.07,0.495){\makebox(0,0)[rt]{$\tilde{H}$\strut}}
  \put(0.99,0.25 ){\makebox(0,0)[rt]{$\tau$\strut}}
\end{picture}
\begin{center}
\begin{minipage}{15cm}
    \caption{Takens-Bogdanov bifurcations without parameters. Case (a)
      hyperbolic; case (b) elliptic. Top: stable and  unstable manifolds;
      bottom: invariant sets. For coordinates and fixed parameters see text.}
    \label{fig2.3}
\end{minipage}
\end{center}
  \end{minipage}}}}
\end{picture}
\end{center}
\end{figure}

Our fourth example addresses Takens-Bogdanov bifurcations without
parameters. In suitable rescaled form it reads
\begin{equation}\begin{array}{rcl}
   \dddot{y} + y\dot{y} &=& \varepsilon ((\lambda -y)\ddot{y}+b \dot{y}^2)
\label{2.4}\end{array}\end{equation}
with fixed parameters $b$, $\lambda$ and $\varepsilon$. The $y$-axis as
equilibrium line is complemented by the two transverse directions
$x=(\ddot{y}, \dot{y})$. Note the algebraically triple zero eigenvalue, double
in the transverse $x$-directions, for $\lambda=y=0$. Two examples of the
resulting dynamics for small positive $\varepsilon$ are summarized in Figure
\ref{fig2.3}.

The coordinates $\tau$ and $\tilde{H}$ in Figure \ref{fig2.3} are adapted to
the completely integrable case $\varepsilon=0$. Indeed, obvious first
integrals are then given by $\Theta =\ddot{y}+\frac{1}{2} y^2$ and
$H=\frac{1}{2} \dot{y}^2 - y \ddot{y} -\frac{1}{3} y^3$. Coordinates are
$\tau=\log \Theta$ and $\tilde{H}=\Theta^{-3/2}H$, not drawn to
scale. Parameters are $\varepsilon, \lambda >0$ and, for the hyperbolic case,
$-17/12 <b<-1$. For the elliptic case we consider $b>-1$. The equilibrium
$y$-axis, a cusp in $(\Theta,H)$ coordinates, transforms to the top (saddles)
and bottom (foci) horizontal boundaries $\tilde{H}=\pm \frac{2}{3}\sqrt{2}$,
with $y=0$ relegated to $\tau=-\infty$. Since $\tau$ and $\tilde{H}$ are
constants of the flow, for $\varepsilon =0$, they represent slow drifts on the
unperturbed periodic motion, for small $\varepsilon >0$ and $\mid
\tilde{H}\mid <\frac{2}{3}\sqrt{2}$. The top value $\tilde{H}
=+\frac{2}{3}\sqrt{2}$ also represents homoclinics to the saddles, for
$\varepsilon =0$.

Along the focus line $\tilde{H}=-\frac{2}{3}\sqrt{2}$ we observe Hopf
bifurcations without parameters, corresponding to $y=\lambda >0$. The value of
$b$ distinguishes elliptic and hyperbolic cases. In addition, lines of
saddle-focus heteroclinic orbits and isolated saddle-saddle heteroclinics are
generated, for $\varepsilon >0$, by breaking the homoclinic sheets of the
integrable case. Note in particular the infinite swarm of saddle-saddle
heteroclinics, in the hyperbolic case.

\begin{figure}
\setlength{\unitlength}{0.88\textwidth}
\begin{center}
\begin{picture}(1.0,0.8)(0,0)
  \put(0.02,0.78){\makebox(0,0)[lt]{(a)}}
  \put(0,0){\makebox(1.0,0.8)
   {
    \includegraphics[width=1.0\unitlength]{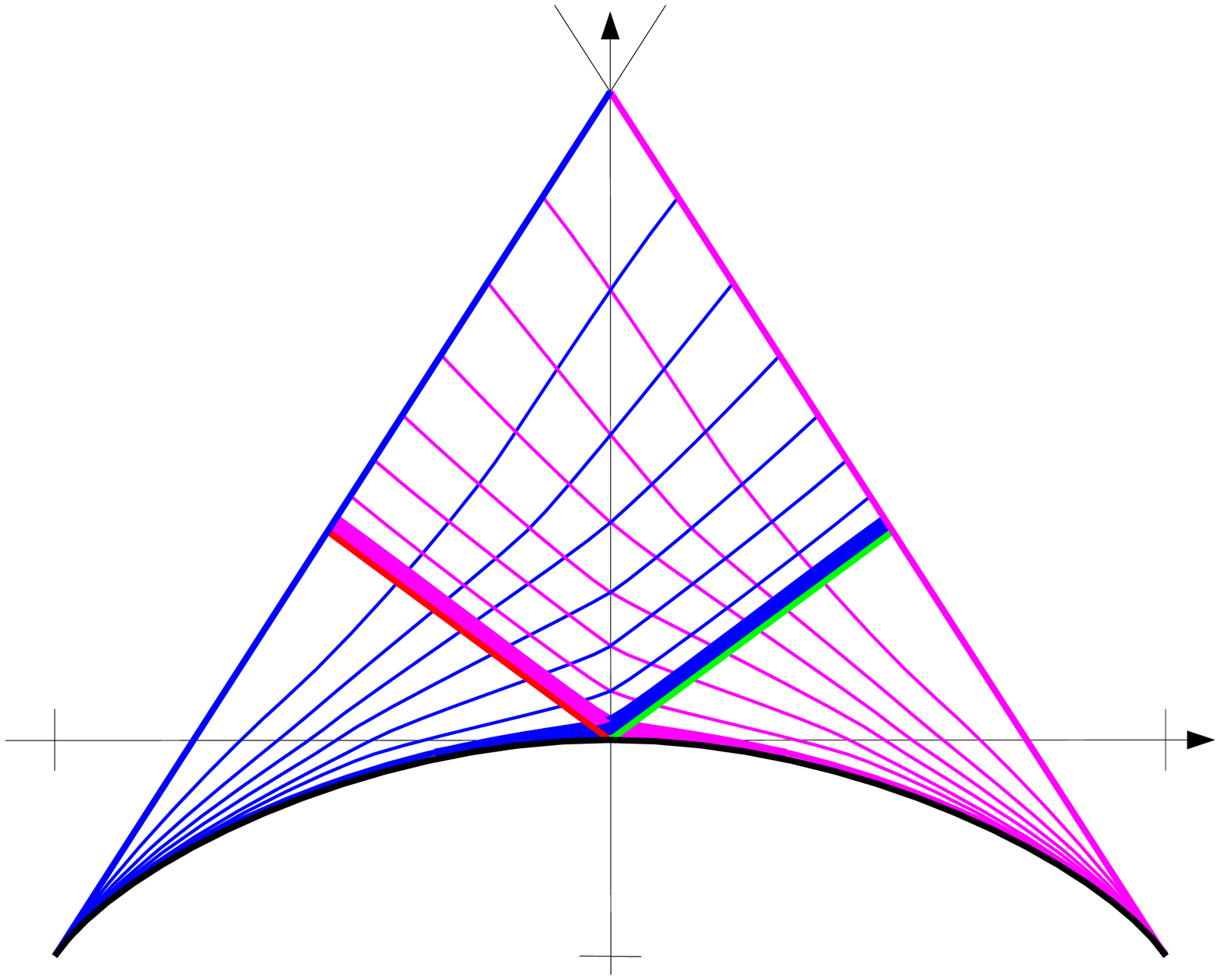}
   }}
  \put(0,0){\framebox(1.0,0.8){~}}
  \put(0.52,0.80){\makebox(0,0)[lt]{$H$\strut}}
  \put(0.99,0.18){\makebox(0,0)[rt]{$\Theta$\strut}}
  \put(0.94,0.23){\makebox(0,0)[b]{$2\sqrt{3}/9$\strut}}
  \put(0.08,0.23){\makebox(0,0)[b]{$-2\sqrt{3}/9$\strut}}
  \put(0.47,0.03){\makebox(0,0)[r]{$-1/12$\strut}}
  \put(0.47,0.72){\makebox(0,0)[r]{$1/4$\strut}}
  \put(0.5,0.18){\makebox(0,0)[t]{hyp.~Hopf\strut}}
  \put(0.5,0.5){\makebox(0,0)[t]{
    \shortstack{infinite\\heteroclinic\strut\\swarm}}}
\end{picture}
\\[0.3ex]
\begin{picture}(1.0,0.8)(0,0)
  \put(0.02,0.78){\makebox(0,0)[lt]{(b)}}
  \put(0,0){\makebox(1.0,0.8)
   {
    \includegraphics[width=1.0\unitlength]{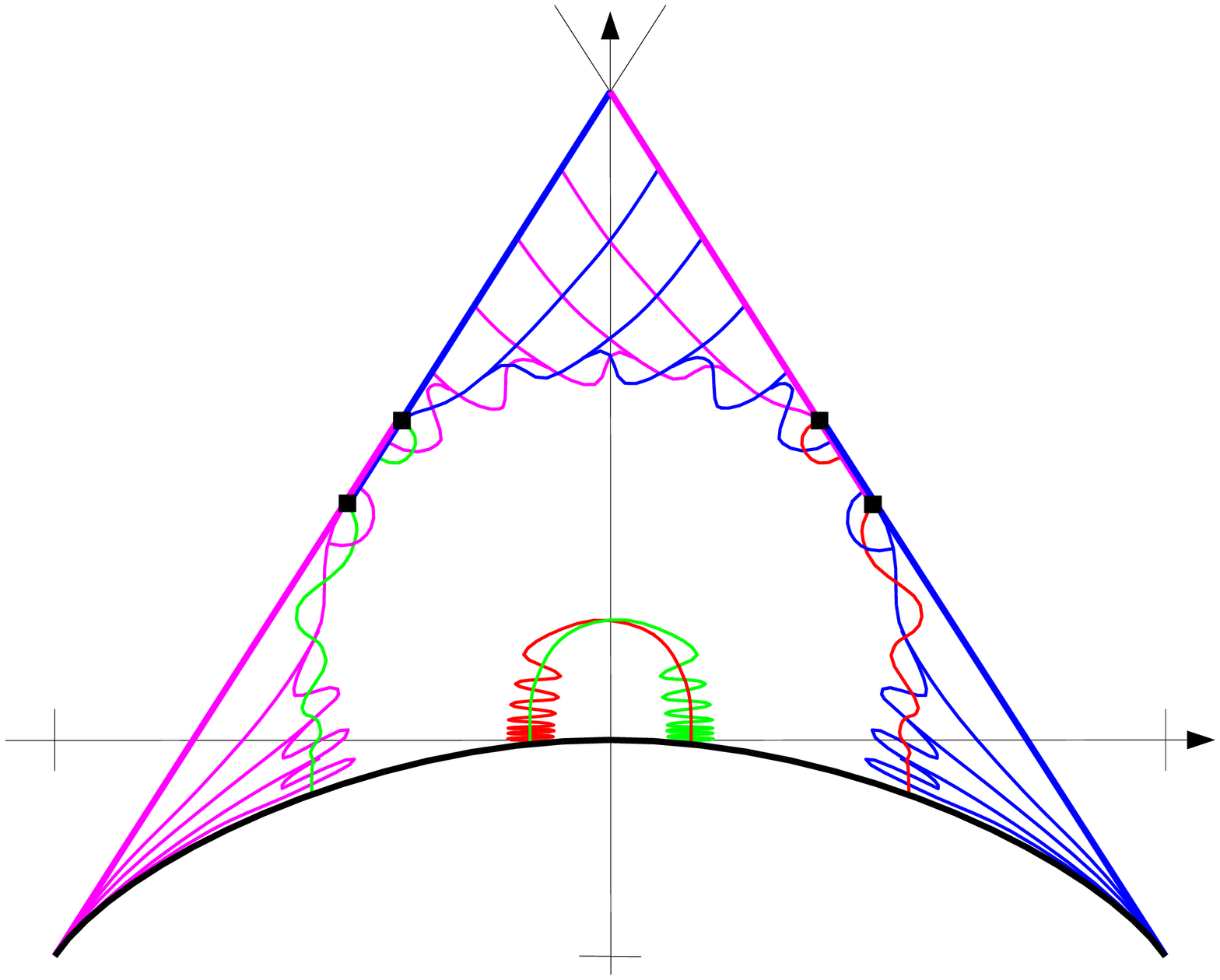}
   }}
  \put(0,0){\framebox(1.0,0.8){~}}
  \put(0.52,0.80){\makebox(0,0)[lt]{$H$\strut}}
  \put(0.99,0.18){\makebox(0,0)[rt]{$\Theta$\strut}}
  \put(0.94,0.23){\makebox(0,0)[b]{$2\sqrt{3}/9$\strut}}
  \put(0.08,0.23){\makebox(0,0)[b]{$-2\sqrt{3}/9$\strut}}
  \put(0.47,0.03){\makebox(0,0)[r]{$-1/12$\strut}}
  \put(0.47,0.72){\makebox(0,0)[r]{$1/4$\strut}}
  \put(0.5,0.18){\makebox(0,0)[t]{ell.~Hopf\strut}}
  \put(0.5,0.63){\makebox(0,0)[t]{
    \shortstack{finite\\heteroclinic\strut\\swarm}}}
\end{picture}
\begin{minipage}{10cm}
\caption{Reversible Takens-Bogdanov bifurcations without parameters. Case (a)
  hyperbolic; case (b) elliptic. For coordinates and fixed parameters see
  text.}
\label{fig2.4}
\end{minipage}
\end{center}
\end{figure}

As a final, fifth example we consider a reversible Takens-Bogdanov bifurcation
without parameters:
\begin{equation}\begin{array}{rcl}
   \dddot{y}+ (1-3y^2)\dot{y} &=& ay \ddot{y}+b \dot{y}^2.
\label{2.5}\end{array}\end{equation}
Here we fix $a,b$ to be small. Again $x=(\ddot{y}, \dot{y})$ denotes the
directions transverse to the equilibrium $y$-axis. Time reversibility
generates solutions $-y(-t)$ from solutions $y(t).$ For two examples of the
resulting dynamics see Figure \ref{fig2.4}. Coordinates are the obvious first integrals, for $\varepsilon =0$, given by $\Theta =\ddot{y}+y-y^3$ and
$H=-\ddot{y} y +\frac{1}{2}\dot{y}^2 + \frac{3}{4} y^4 - \frac{1}{2}
y^2$. Note the two Takens-Bogdanov cusps, separated by a Hopf point along the
lower arc of the equilibrium ``triangle''. Compare Figures \ref{fig2.2},
\ref{fig2.3}. The elliptic Hopf point (b) arises for $a\cdot(a-b)>0$, whereas
$a\cdot(a-b)<0$ in the hyperbolic case (a). Also note the associated finite
and infinite swarms of saddle-saddle heteroclinics, respectively.

\Section{Methods} \label{sec3} \setzero \addsec

\vskip-5mm \hspace{5mm}

Pictures are not proofs. What has been proved, then, and how? We use
ingredients involving algebra, analysis, and numerical analysis, as we outline
in this section. For further details see \cite{Liebscher00-Dissertation,
  FiedlerLiebscherAlexander98-HopfTheory, FiedlerLiebscher01-TakensBogdanov,
  AfendikovFiedlerLiebscher02-PlaneKolmogorovFlows}.

In a first algebraic step, we derive normal forms for vector fields with lines
or planes of equilibria, assuming spectral degeneracies of the linearization
$A$ in transverse directions. The spectral assumptions on $A$ in fact coincide
with those established for parameter-dependent matrix families in standard
bifurcation theory. This is reflected in the naming of the five examples of
section \ref{sec2}.

There are more or less standard procedures to derive normal forms of vector
fields. By suitable polynomial diffeomorphisms, certain Taylor coefficients of
the vector field are eliminated, successively, for higher and higher
order. See for example \cite{Vanderbowhede89-CentManifolds} for a systematic
choice of normal forms, particularly apt for introducing equivariance of the
nonlinear normal-form terms under the action of $\exp(A^Tt).$ Normal forms
are, however, nonunique in general and other choices are possible.

In the present cases, we adapt the normal-form procedure to preserve the
locally flattened equilibrium manifolds. Although the approach in
\cite{Vanderbowhede89-CentManifolds} can be modified to accommodate that
requirement, it did not provide vector fields convenient for subsequent
analysis of the flow. A systematic approach to this combined problem is not
known, at present. All examples (\ref{2.1})--(\ref{2.5}) represent truncated
normal forms. For the derivation of specific normal forms, for example of
(\ref{2.4}), to any order, and of example (\ref{2.5}), to third order, see
\cite{FiedlerLiebscher01-TakensBogdanov},
\cite{AfendikovFiedlerLiebscher02-PlaneKolmogorovFlows}, respectively.

Subsequent analysis of the normal-form vector fields is based on scalings,
alias blow-up constructions. This is the origin, for example, of the small
scaling parameter $\varepsilon$ in the Takens-Bogdanov example (\ref{2.4}). In
passing, we note a curious coincidence of two view points for (\ref{2.4}),
concerning the roles of the equilibrium coordinate $y\in\setR$ and the
``fixed'' real parameter $\lambda \in\setR$. First, we may consider $\lambda$
as a parameter, with a line of equilibria $y$ associated to each fixed
$\lambda$. Then (\ref{2.4}) describes the collision of a transverse zero
eigenvalue at $y=0$, as in (\ref{2.1}), with imaginary Hopf eigenvalues at
$y=\lambda >0$, as in (\ref{2.3}), as $\lambda$ decreases through
zero. Alternatively, we may consider normal forms for a plane
$\mathbf{y}=(y_1,y_2)$ of equilibria with a transverse double zero eigenvalue,
at $\mathbf{y}=0$. It turns out that the two cases coincide, after a scaling
blow-up, up to second order in $\varepsilon$, via the
correspondence $y=y_1,\lambda =y_2$.

The core of any successful flow analysis in bifurcation theory is an
integrable vector field; see again section \ref{sec2}. The issues of
nonintegrable perturbations, by small $\varepsilon>0$, and of omitted higher
order terms, not in normal form, both ensue. Since the underlying integrable
dynamics is periodic or homoclinic, in examples (\ref{2.3})--(\ref{2.5}),
averaging procedures apply. Indeed, $\varepsilon >0$ then introduces a
periodically forced, slow flow on first integrals, like $(\Theta,H)$,
characteristic of $\varepsilon=0$. We therefore derive an appropriate, but
autonomous Poincar\'e flows, on $(\Theta,H)$, such that the associated true
Poincar\'e map can be viewed as a time discretization of first order and step
size $\varepsilon$. In the unperturbed periodic region, this amounts to
averaging, while the Poincar\'e flow indicates Melnikov functions at
homoclinic or heteroclinic boundaries. The exponential averaging results by
Neishtadt \cite{Neishtadt84-ExpSmalSplit}, for example, then imply that the
separatrix splittings in the elliptic Hopf case (b), indicated in Figure
\ref{fig2.2}, are exponentially small in the radius $r$ of the split sphere,
for analytic vector fields. See also \cite{FiedlerScheurle96-HomoclinicOrbits}.

Lower bounds of separatrix splittings have not been established, in our
settings. This problem is related to the very demanding Lazutkin program of
asymptotic expansions for exponentially small separatrix splittings. For
recent progress, including the case of Takens-Bogdanov bifurcations for
analytic maps, see \cite{Gelfreich99-ExpSeparatices} and the references
there. In absence of rigorous lower bounds, our figures indicate only simplest
possible splitting scenarios.

While the splitting near elliptic Hopf points are exponentially small, the
discretization of the Poincar\'e  flow also exhibits splittings of the
unperturbed saddle homoclinic families, which are of first order in the
perturbation parameter $\varepsilon $, in example (\ref{2.4}), or in the small
parameters $a, b$, in example (\ref{2.5}). Explicit expressions have been
derived for the Melnikov functions associated to these homoclinic splittings,
in terms of elliptic function in case (\ref{2.4}), and even of elementary
functions in case (\ref{2.5}). Simplicity and uniqueness of zeros of the
Melnikov functions, however, has only been confirmed numerically. While this
does not, strictly speaking, match an analytic proof, it still at least
supports the validity of the scenarios summarized in Figures \ref{fig2.3} and
\ref{fig2.4}.

\Section{Interpretation and perspective} \label{sec4} \setzero \addsec

\vskip-5mm \hspace{5mm}

We indicate some consequences of the above results for the examples of coupled
oscillators, viscous shock profiles, and Kolmogorov flows indicated in section
\ref{sec1}. We conclude with a few remarks on the future perspective of
bifurcations without parameters.

We first return to the example (\ref{1.4})--(\ref{1.8}) of a coupled
oscillator square, $m=1.$ Equilibria $y$ of the Poincar\'e flow then indicate
decoupled antipodal periodic pairs, say with phase difference $y+c.$ The case
of a transverse zero eigenvalue, (\ref{2.1}) and Figure \ref{fig2.1}, then
indicates a 50\% chance of recovery of decoupling with a stable phase
difference $y<0$, locally, even when the stability threshold $y=0$ has been
exceeded. The hyperbolic Hopf case  (\ref{2.2}), Figure \ref{fig2.2} (a),
illustrates immediate oscillatory loss of decoupling stability by transverse
imaginary eigenvalues. A 100\% recovery of decoupling stability, in contrast,
occurs at elliptic Hopf points; see (\ref{2.2}), Figure \ref{fig2.2} (b). The
exponentially small Neishtadt splitting of separatrices indicates a very
delicate variability in the asymptotic phase relations of this recovery, for
$t\rightarrow\pm\infty$. See \cite{FiedlerLiebscherAlexander98-HopfTheory}.
The Takens-Bogdanov cases (\ref{2.4}), Figures \ref{fig2.3} (a), (b) can then
be viewed as consequences of a mutual interaction, of a transverse zero
eigenvalue with either Hopf case, for recovery of stable decoupling.

In the example (\ref{1.9}), (\ref{1.10}) of stiff balance laws, elliptic Hopf
bifurcation without parameters as in (\ref{2.4}), Figure \ref{fig2.3} (b),
indicates {\em oscillatory} shock profiles $u(\tau), \tau=(\xi
-st)/\varepsilon.$ Such profiles in fact contradict the Lax condition, being
over-compressive, and violate standard monotonicity criteria. For small
viscosities $\varepsilon >0$, weak viscous shocks in fact turn out unstable,
in any exponentially weighted norm, unless they travel at speeds $s$ exceeding
all characteristic speeds. The oscillatory profiles can be generated, in fact,
by the interaction of inherently non-oscillatory gradient flux functions
$F(u)$ with inherently non-oscillatory gradient-like kinetics $G(u)$ in
systems of $\dim \geq 3$. See \cite{Liebscher00-Dissertation}.

The problem of plane stationary Kolmogorov flows asks for stationary solutions
of the incompressible Navier-Stokes equations in a strip domain $(\zeta,
\eta)\in\setR\times [0,2\pi]$, under periodic boundary conditions in $\eta $;
see \cite{MeshalkinSinai61-PlaneFluidFlow}.
An $\eta$-periodic external force $(F(\eta),0)$, is
imposed, acting in the unbounded $\zeta $-direction. Kolmogorov chose
$F(\eta)=\sin \eta$. The Kirchg\"assner reduction
\cite{Kirchgaessner82-Reduction} captures all
bounded solutions which are nearly homogeneous in $\zeta$, in a center
manifold spirit which lets us interpret $\zeta$ as ``time''. The resulting
ordinary differential equations in $\setR^6$ reduce to $\setR^3$, by fixing
the values of three first integrals. A line of $\zeta$-homogeneous equilibria
appears, in fact, and Kolmogorov's choice corresponds to example (\ref{2.5})
with  $a=b=0$, to leading orders. In particular note the double reversibility,
then, under $y(t)\mapsto \pm y(-t)$ which is generated by
\begin{equation}\begin{array}{rcl}
   F(\eta) &=& -F(\eta +\pi), \qquad \mbox{and}
\\ F(\eta) &=& -F(-\eta).
\label{4.1} \end{array}\end{equation}
As observed by Kolmogorov, an abundance of spatially periodic profiles results.
The sample choice $F(\eta)= \sin \eta + c \sin 2\eta$, in contrast, which
breaks the first of the symmetries in (\ref{4.1}), leads to (\ref{2.5}) with
$b=0<a$, alias an elliptic reversible Takens-Bogdanov point without
parameters; see Figure \ref{fig2.4} (b). In particular, the set of
near-homogeneous bounded velocity profiles of the incompressible, stationary
Navier-Stokes system is then characterized by an abundance of oscillatory
heteroclinic wave fronts, which decay to different asymptotically homogeneous
$\zeta$-profiles, for $\zeta \rightarrow \pm \infty$. The PDE stability of
these heteroclinic profiles is of course wide open.

As for perspectives of our approach, we believe to have examples at hand, from
sufficiently diverse origin, to justify further development of a theory of
bifurcations without parameters. In fact, transverse spectra $\{0, \pm i
\omega\}$ and $\{\pm i \omega_1, \pm i \omega_2\}$ still await investigation
before we can claim any insight into nonhyperbolicity of even the simple case
of an equilibrium plane. This assumes the absence of further structural
ingredients like symplecticity, contact structures, symmetries, and the
like. Certainly our example collecting activities are far from complete, at
this stage.

In addition, we have not addressed the issue of perturbations, so far, which
could destroy the equilibrium manifolds by small drift terms. Examples arise,
for example, when slightly detuning the basic frequencies $2\pi/T_j$ of our
uncoupled oscillators or, much more generally, in the context of multiple
scale singular perturbation problems. Feedback and input from our readers will
certainly be most appreciated!


\label{lastpage}

\end{document}